\documentclass[10pt]{amsart}
\usepackage[pagebackref]{hyperref}

\renewcommand{\d}{\mathrm{d}}

\newcommand{\ts}{\textstyle }\newcommand{\ds}{\displaystyle }

\newcommand{\cI}{{\mathcal I}}
\newcommand{\bbR}{{\mathbb R}}
\newcommand{\bbC}{{\mathbb C}}

\newcommand{\bbT}{{\mathbb T}}
\newcommand{\bbZ}{{\mathbb Z}}

\newcommand{\E}{{\mathrm e}}
\newcommand{\iC}{{\mathrm i}}

\newcommand{\SO}{\operatorname{SO}}
\newcommand{\SU}{\operatorname{SU}}
\newcommand\Un{{\operatorname U}}

\newcommand{\Ad}{\operatorname{Ad}}

\newcommand{\p}{{\partial}}

\DeclareMathOperator{\tr}{tr}

\newcommand{\eug}{{\mathfrak{g}}}
\newcommand{\eusu}{{\mathfrak{su}}}
\newcommand{\euso}{{\mathfrak{so}}}

\newcommand{\euz}{{\mathfrak{z}}}

\newcommand{\la}{\langle}
\newcommand{\ra}{\rangle}

\newcommand{\w}{{\mathchoice{\,{\scriptstyle\wedge}\,}{{\scriptstyle\wedge}}
      {{\scriptscriptstyle\wedge}}{{\scriptscriptstyle\wedge}}}}

\numberwithin{equation}{section}

\newtheorem{theorem}{Theorem}

\newtheorem{proposition}{Proposition}

\theoremstyle{remark}

\newtheorem{remark}{Remark}
\newtheorem{example}{Example}

\begin{document}

\author[R. Bryant]{Robert L. Bryant}
\address{Duke University Mathematics Department\\
         P.O. Box 90320\\
         Durham, NC 27708-0320}
\email{\href{mailto:bryant@math.duke.edu}{bryant@math.duke.edu}}
\urladdr{\href{http://www.math.duke.edu/~bryant}%
         {http://www.math.duke.edu/\lower3pt\hbox{\symbol{'176}}bryant}}

\title[$SO(n)$-invariant special Lagrangians with fixed loci]
      {$\SO(n)$-invariant Special Lagrangian\\
       Submanifolds of $\bbC^{n+1}$\\
       with Fixed Loci}

\date{March 25, 2004}

\begin{abstract}
Let~$\SO(n)$ act in the standard way on~$\bbC^n$
and extend this action in the usual way to~$\bbC^{n+1} 
= \bbC\oplus\bbC^n$.  

It is shown that a nonsingular special Lagrangian
submanifold~$L\subset\bbC^{n+1}$ that is invariant
under this~$\SO(n)$-action intersects the fixed~$\bbC
\subset\bbC^{n+1}$ in a nonsingular real-analytic 
arc~$A$ (which may be empty).  
If~$n>2$, then $A$ has no compact component.

Conversely, an embedded, noncompact nonsingular 
real-analytic arc~$A\subset\bbC$ lies in an embedded
nonsingular special Lagrangian submanifold that is 
$\SO(n)$-invariant. The same existence result 
holds for compact~$A$ if~$n=2$.  If~$A$ is connected, 
there exist $n$ distinct nonsingular $\SO(n)$-invariant
special Lagrangian extensions of~$A$ such that
any embedded nonsingular $\SO(n)$-invariant
special Lagrangian extension of~$A$ agrees with
one of these~$n$ extensions in some open neighborhood
of~$A$.  

The method employed is an analysis of a singular nonlinear
\textsc{pde} and ultimately calls on the work of G\'erard
and Tahara to prove the existence of the extension.
\end{abstract}



\thanks{
Thanks to Duke University for its support via
a research grant, to the NSF for its support
via DMS-0103884, to the Mathematical Sciences
Research Institute, and to Columbia University.  
\hfill\break
\hspace*{\parindent} 
This is Version~$2$.  The first version was posted to the
arXiv on 12 February 2004.
}

\maketitle

\setcounter{tocdepth}{2}
\tableofcontents

\section{Introduction}\label{sec: intro}

\subsection{Special Lagrangian geometry}\label{ssec: intro sLag}
In the now-classic 1982 paper~\cite{MR1985i:53058} of Harvey
and Lawson, an $m$-dimensional submanifold~$L\subset\bbC^m$ 
is defined to be \emph{special Lagrangian} if it is
Lagrangian with respect to the standard K\"ahler form
\begin{equation}
\omega = \frac{\iC}2\bigl(\d z_1\w \d \overline{z_1}
+ \cdots + \d z_m\w \d \overline{z_m}\bigr)
\end{equation}
and, moreover, the standard holomorphic volume 
form~$\d z = \d z_1 \w \cdots \d z_m$ pulls back to~$L$ 
to be a real $n$-form, i.e., $\Upsilon =\text{Im}(\d z)$ pulls back
to~$L$ to be identically zero.%
\footnote{In fact, Harvey and Lawson show that
if~$L\subset\bbC^m$ is Lagrangian and oriented, then the
pullback of~$\d z$ to~$L$ is of the form~$\E^{\iC\phi}\,\text{vol}_L$,
for some function~$\phi:L\to S^1$, where~$\text{vol}_L$ is the
$n$-form on~$L$ that is the volume form of the induced metric
on~$L$.  The `phase factor'~$\E^{\iC\phi}$ is a sort of 
complex determinant, so the equation $\E^{\iC\phi}=1$
can be thought of as setting a determinant equal to~$1$, 
hence the modifier `special'.}

In other words, an $m$-dimensional submanifold~$L\subset\bbC^m$
is special Lagrangian if and only if it is an integral manifold
of the ideal~$\cI$ generated by the differential forms~$\omega$
and~$\Upsilon$.

As Harvey and Lawson show, an $L\subset\bbC^m$ is
special Lagrangian if and only if it is calibrated by
the $m$-form~$\Phi = \text{Re}(\d z)$.  In particular,
special Lagrangian submanifolds are mass-minimizing.

In the intervening 20 years, special Lagrangian submanifolds and
foliations whose leaves are special Lagrangian submanifolds
(and their generalizations to Calabi-Yau manifolds in 
the place of~$\bbC^m$) have turned out to be important in
several areas in differential geometry and theoretical physics.
The interested reader can consult recent papers of Joyce,
such as~\cite{math.DG/0111111} for a survey and references
to the (by now rather large) literature on this subject.

\subsection{Symmetry reduction}\label{ssec: intro symred}
Because understanding the possible types of singularities
of special Lagrangian submanifolds is important for applications,
considerable effort has been expended to construct
explicit examples.

One of the standard methods of constructing explicit integral
manifolds of an ideal~$\cI$ is to use symmetry reduction, i.e., 
to fix a subgroup~$G$ of the symmetries of the given ideal~$\cI$
and look for solutions that are invariant under the action of~$G$.

\subsubsection{Symmetries of~$\cI$}\label{sssec: idealsyms}
In the case of special Lagrangian geometry, the symmetry
group of the ideal~$\cI$ depends to some extent 
on the dimension~$m$.  

The dimension~$m=2$ is exceptional.  In this case, there exists
a complex structure~$J$ on~$\bbC^2$ (different from the standard one)
such that~$\omega$ and~$\Upsilon$ are the real and imaginary 
parts, respectively, of a $J$-holomorphic $(2,0)$-form on~$\bbC$.
Thus, the local symmetries of~$\cI$ are simply the local
biholomorphisms of~$(\bbC^2,J)$.  (Moreover, the special Lagrangian
submanifolds are simply the $J$-complex curves in~$(\bbC^2,J)$.)

For~$m>2$, the symmetry group of the ideal~$\cI$ is 
finite dimensional and consists of the group generated by 
the translations in~$\bbC^m$, 
the group~$\SU(m)$ acting linearly on~$\bbC^m$, 
the dilations by nonzero complex numbers~$\lambda$ 
that satisfy~$\lambda^m\in\bbR$, 
and conjugation.    

\subsubsection{Cohomogeneity one examples}\label{sssec: cohom1ex}
Harvey and Lawson themselves considered and solved the problem
of describing the~$G$-invariant special Lagrangian submanifolds
of~$\bbC^m$ when~$G\subset\SU(m)$ is either~$\SO(m)$ 
or~$\bbT^{m-1}$, a maximal torus in~$\SU(m)$.  
The Lagrangian-isotropic orbits of these actions have 
dimension~$m{-}1$, and so the standard theory leads one to 
expect that, in these cases, the $G$-invariant special 
Lagrangian manifolds will be found by solving a single 
\textsc{ode}.  

Indeed, this is precisely what Harvey and Lawson 
find.  They show that, when~$G=\SO(m)$, the sets
\begin{equation}
L_c = \{ \zeta\,{\bf u}\ \vrule\
 {\bf u}\in S^{m-1},\, \zeta\in\bbC,\ \text{Im}(\zeta^n) = c\}
\end{equation} 
for~$c$ a real constant are $\SO(m)$-invariant and
special Lagrangian at smooth points. The set~$L_0$
is the union of $m$ special Lagrangian $m$-planes
that are each individually~$\SO(m)$-invariant.
When~$c\not=0$, the set~$L_c$, which is smooth, 
is the disjoint union of~$m$ connected components,
each of which is diffeomorphic to~$\bbR\times S^{m-1}$ 
with each `end' asymptotic to a $SO(m)$-invariant 
special Lagrangian $m$-plane.  Moreover, any smooth, 
connected $\SO(m)$-invariant special Lagrangian submanifold 
of~$\bbC^m$ is an open subset of one of these examples.

In the case of~$G=\bbT^{m-1}$, the picture is slightly 
more complicated:  The $\bbT^{m-1}$-invariant special 
Lagrangian submanifolds are the simultaneous level sets
of the functions~$f_0,\ldots,f_{m-1}$ where
\begin{equation}
\begin{aligned}
f_0 &= \text{Im}\bigl(z_1z_2\cdots z_m\bigr)\\
f_k &= |z_k|^2 - |z_m|^2, \qquad 1 \le k < m.
\end{aligned}
\end{equation}
Many of these level sets are singular (some are
analytically irreducible and some are not) and they 
furnish interesting examples of the kinds of singularities
that mass minimizing currents can display.  (In fact, 
this is one of the reasons that Harvey and Lawson found
them so interesting.)

\subsubsection{General compact group actions}
\label{sssec: genG}
It may clarify matters to consider special 
Lagrangian symmetry reduction in the case of an 
arbitrary connected compact subgroup~$G$, that,
without loss of generality, can be assumed to be 
a subgroup of~$\SU(m)$.

Recall that the Lie algebra $\eusu(m)$ of~$\SU(m)$ 
is the vector space of $m$-by-$m$ skew-Hermitian 
matrices and is endowed with a positive definite 
$\Ad\bigl(\SU(m)\bigr)$-invariant inner product 
defined by~$\la a,b\ra = -\tr\bigl(ab\bigr)$.

Let~$\eug\subset\eusu(m)$ denote the
Lie algebra of~$G$ and write~$\eug=\euz\oplus[\eug,\eug]$
where~$\euz\subset\eug$ is the center of~$\eug$. 
Using the above inner product,~$\eug$ can be
identified with its dual, so~$\eug^*$ and $\eug$ 
will be identified henceforth.  
Let~${\pi_\eug}:\eusu(m)\to\eug$ denote the orthogonal
projection onto~$\eug$.  

The linear action of~$G$ on~$\bbC^m$ is $\omega$-Poisson, 
with momentum mapping~$\mu:\bbC^m\to\eug=\eug^*$ 
given by
\begin{equation}
\mu(z) = {\pi_\eug}\bigl(\,\iC\,z\,{}^t\bar z\,\bigr).
\end{equation}
Note that a~$G$-orbit~$G\cdot z\subset\bbC^m$ is
$\omega$-isotropic if and only if $\mu(z)$ 
lies in~$\euz$.

\subsubsection{Smooth reduction}\label{sssec: smthred}
For a fixed~$\xi\in\euz$, 
let~$\mu^{-1}(\xi)^*\subset\mu^{-1}(\xi)\subset\bbC^m$ 
denote the subset that consists of $\mu$-clean points.%
\footnote{For a smooth mapping~$f:X\to Y$
between smooth manifolds, a point~$x\in X$ is said
to be~\emph{$f$-clean} 
if the level set~$Z = f^{-1}\bigl(f(x)\bigr)\subset X$ 
is a smooth submanifold near~$x$ 
and~$T_xZ=\ker f'(x)$.}

The submanifold~$\mu^{-1}(\xi)^*$ is usually dense 
in~$\mu^{-1}(\xi)$ and has dimension~$m{+}k$ 
for some~$k\ge 1$. The $G$-orbits in~$\mu^{-1}(\xi)^*$ 
are~$\omega$-isotropic and of dimension~$m{-}k$.  

The symplectic quotient~$M^*_\xi=G\backslash\mu^{-1}(\xi)^*$
will, at most points, be a smooth symplectic manifold
of dimension~$2k$, with the canonical projection
$\pi_\xi:\mu^{-1}(\xi)^*\to M^*_\xi$ being
a smooth submersion and the induced symplectic
form~$\omega_\xi$ on~$M^*_\xi$ having the property
that~$\pi^*_\xi(\omega_\xi)$ is equal to the
pullback of~$\omega$ to~$\mu^{-1}(\xi)^*$.

The reduced space $M^*_\xi$ inherits a metric~$g_\xi$ 
defined by the condition that~$\pi_\xi:\mu^{-1}(\xi)^*\to M^*_\xi$ 
be a Riemannian submersion.  The pair~$(g_\xi,\omega_\xi)$
then define a K\"ahler structure on~$M^*_\xi$.

For any $\omega_\xi$-Lagrangian submanifold~$L_\xi\subset M_\xi$,
its preimage~${\pi_\xi}^{-1}(L_\xi)\subset\mu^{-1}(\xi)^*\subset\bbC^m$
is a $G$-invariant $\omega$-Lagrangian submanifold of~$\bbC^m$.  

Moreover, it is not difficult to show that there is 
a~$k$-form~$\Upsilon_\xi$ defined on~$M^*_\xi$ (only
up to a sign if~$\mu^{-1}(\xi)^*$ is not orientable)
with the property that~$\Upsilon_\xi$ vanishes 
when pulled back to~$L_\xi$ if and only if
${\pi_\xi}^{-1}(L_\xi)$ is special Lagrangian.
In fact, up to an orientation,~$\Upsilon_\xi$ can
be seen as the imaginary part of a $(k,0)$-form
on the K\"ahler manifold~$M^*_\xi$. 

Thus, away from the singularities of the various 
mappings, the problem of describing the $G$-invariant
special Lagrangian submanifolds of~$\bbC^m$ can be
reduced to a similar problem in lower dimensions.  

\begin{example}[Cohomogeneity~$1$]\label{ex: cohom1}
For example, when~$k=1$ (as is true in the Harvey-Lawson
examples), $M^*_\xi$ is a surface and~$\Upsilon_\xi$
is a $1$-form on~$M^*_\xi$. The problem of describing
the $G$-invariant special Lagrangian submanifolds 
(away from singularities) is thus reduced to finding 
the integral curves of a $1$-form on a surface.

As another example, to illustrate the method and because
this will be used to construct a needed example later,
consider the action 
of~$G=S^1\times\SO(p)\times\SO(q)$ on~$\bbC^{p+q}$
acting on~$\bbC^{p+q}$ via the action
\begin{equation}
(e^{\iC\theta},A,B)\cdot(z,w) 
= \bigl(e^{q\iC\theta}Az,\, e^{-p\iC\theta}Bw\bigr),
\qquad\qquad z\in\bbC^p,\,w\in\bbC^q.
\end{equation}
The momentum mapping is
\begin{equation}
\mu(z,w) = \bigl(\iC(q|z|^2-p|w|^2),
 \text{Im}(z\,{}^t\overline z), \text{Im}(w\,{}^t\overline w)\bigr).
\end{equation}
Taking~$\xi=(\iC c,0,0)$ for some constant~$c\in\bbR$, 
(which is the only allowable choice unless~$p$ or~$q$ is~$2$) yields
\begin{equation}
\mu^{-1}(\xi)
= \left\{\,(\,\zeta\,{\bf u},\,\eta\,{\bf v}\,)
          \,\vrule\ 
       {\bf u}\in S^{p-1},\,{\bf v}\in S^{q-1},\,\zeta,\eta\in\bbC,\,
        q|\zeta|^2-p|\eta|^2=c\,\right\}.
\end{equation}
One then finds that~$\Upsilon_\xi = \d\bigl(\text{Re}(\zeta^q\eta^p)\bigr)$,
so that~$\mu^{-1}(\xi)$ is `foliated' (some of the `leaves' may be 
singular) by the level sets of the function~$\text{Re}(\zeta^q\eta^p)$,
which are special Lagrangian.
\end{example}

\begin{example}[Cohomogeneity~$2$ and almost complex surfaces]
\label{ex: cohom2}
The $k=2$ case is somewhat more interesting.
In this case, assuming that~$\mu^{-1}(\xi)^*$ is orientable,%
\footnote{The nonorientable case does occur, as will
be seen below.} it is not difficult to show
that~$M^*_\xi$ (which has real dimension~$2k=4$), inherits 
a natural almost complex structure~$J_\xi$ such that~$L_\xi
\subset M^*_\xi$ is a $J_\xi$-complex curve if and 
only if ${\pi_\xi}^{-1}(L_\xi)$ is special Lagrangian.
(In fact, $\omega_\xi$ and~$\Upsilon_\xi$ in this
case turn out to be essentially the real and imaginary
parts of a $(2,0)$-form with respect to the almost
complex structure~$J_\xi$.)  If~$\mu^{-1}(\xi)^*$ is not
orientable, then by passing to its orientation double
cover, one can define~$\Upsilon_\xi$ and~$J_\xi$
as before, but on a covering space of~$M^*_\xi$.
\end{example}

\subsubsection{Singularities in reduction}\label{sssec: singred}
Thus, one understands, in a general way, how to describe
the special Lagrangian submanifolds invariant under a 
compact Lie group, at least away from singular points
of the quotient~$M_\xi = G\backslash \mu^{-1}(\xi)$.
However, in the case of cohomogeneity greater than~$1$,
as in the recent work of Dominic Joyce~\cite{math.DG/0111111}
(cf., especially, the series 
of papers~\cite{math.DG/0111324,math.DG/0111326,math.DG/0204343}),
the singular locus plays an important role.  The standard
approach described above is not adequate to address 
existence and uniqueness questions in this situation.

To take a specific example, Joyce considers the case 
of~$G=\SO(2)$ acting on~$\bbC^3 = \bbC^1\oplus\bbC^2$ 
(trivially on the first summand~$\bbC^1$ and 
in the standard way on the second summand~$\bbC^2$).%
\footnote{The $\SO(2)$-action that Joyce actually considers 
is conjugate to this one in~$\SU(3)$, but this obviously 
will not affect the results.}  

Joyce showed that a nonsingular $\SO(2)$-invariant 
special Lagrangian~$3$-fold~$L\subset\bbC^3$
meets the fixed factor~$\bbC^1$, if at all, 
in a nonsingular real-analytic arc~$A\subset\bbC^1$. 

Now, in fact, $A$ lies in two distinct~$\SO(2)$-invariant, 
special Lagrangian $3$-folds, namely~$L$ and~$\iC{\star}L$, 
where~$\iC$ acts on~$\bbC^3$ as~$\iC{\star}(z_0,z_1,z_2)
= (\,z_0,\,\iC\,z_1, \,\iC\,z_2\,)$.  (That~$\iC{\star}L$
is special Lagrangian and distinct from~$L$ is left 
to the reader, but see~\S\ref{sssec: lambdaact}.)

For use in his work on singularities 
of special Lagrangian $3$-folds, Joyce asked%
\footnote{private communication, 2 January 2002.} 
whether, any nonsingular $\SO(2)$-invariant
special Lagrangian~$3$-fold that meets~$\bbC$ in
the arc~$A\subset\bbC^1$ is equal to one of~$L$
or~$\iC{\star}L$ in some neighborhood of~$A$.

\subsection{Results}\label{ssec: results}
As will be proved in this article, a stronger statement 
is true:  For any nonsingular, connected, embedded, 
real-analytic arc~$A\subset\bbC$,
there exist two distinct $\SO(2)$-invariant, 
nonsingular, embedded, connected special Lagrangian
submanifolds of~$\bbC^3$ that intersect the 
fixed~$\bbC^1$ in the analytic arc~$A$. 
Moreover, any $\SO(2)$-invariant, nonsingular, 
embedded, connected special Lagrangian
submanifold of~$\bbC^3$ that intersects the 
fixed~$\bbC^1$ in the analytic arc~$A$ agrees with
one of these two special Lagrangian submanifolds in
some open neighborhood of~$A$.  

Moreover, these existence and uniqueness
results generalize when appropriately stated, 
when~$A$ is noncompact,%
\footnote{When~$n>2$, the assumption of 
noncompactness is necessary, 
cf. Proposition~\ref{prop: Anonsing}.}
to the case of~$\SO(n)$ 
acting on~$\bbC^{n+1}=\bbC^1\oplus\bbC^n$ (trivially 
on the first summand~$\bbC^1$ and in the standard way 
on the second summand~$\bbC^n$).  As this is no more 
difficult than the case~$n=2$ (which is the one of 
interest to Joyce), these more general results will 
be proved in this article, 
cf. Theorems~\ref{thm: nuniqueL} and~\ref{thm: existL}.

The method used to establish this uniqueness result
is basically an elementary, though slightly delicate, 
examination of a power series expansion.  However, the 
this argument, while it proves uniqueness, 
only suffices to prove existence in the category 
of formal power series.  
It is not adequate to address the problem of existence 
and the elementary argument does not seem to lend itself
to any of the usual methods of proving convergence of
the formal power series.  However, it turns out that,
by calling on the work of G\'erard and Tahara on existence
of solutions of certain kinds of singular holomorphic 
\textsc{pde}, the convergence of this power series 
can be established.

\subsection{Acknowledgements}\label{ssec: acknowl}
This note resulted from questions that Dominic Joyce asked
about group-invariant special Lagrangian $3$-folds.
It is a pleasure to thank him for raising these questions
and for several discussions about their significance
as well as for comments on an early draft of this article.
The key local existence argument making use of the work 
of G\'erard and Tahara was made possible by Professor Tahara's 
generous help and advice, which is gratefully acknowledged.

After Version~1 of this article was posted to the arXiv
on~12 February 2004, Professors Castro and Urbano made me aware 
of their preprint~\cite{math.DG/0403108}, in which nonsingular 
reduction in the case of~$\SO(p)\times\SO(q)$ acting on~$\bbC^{p+q}$ 
is discussed at greater length.  I thank them for bringing
their article to my attention.

\section{Invariant special Lagrangians}
\label{sec: invsLags}

Write~$z_j = x_j + \iC\,y_j$ ($0\le j\le n$) for the real and imaginary 
parts of the complex coordinates in~$\bbC^{n+1}$.  The K\"ahler form is 
\begin{equation}
\omega = {\ts\frac{\iC}2}\,\bigl(\d z_0\w \d \overline{z_0}
+ \d z_1\w \d \overline{z_1}
+ \cdots + \d z_n\w \d \overline{z_n}\bigr)
\end{equation}
and the imaginary part of the complex volume form is
\begin{equation}
\Upsilon = \text{Im}\bigl(\d z_0 \w\d z_1 \w \cdots \d z_m\bigr)
\end{equation}

By definition, a (real) $(n{+}1)$-dimensional 
submanifold~$L\subset\bbC^{n+1}$ is special Lagrangian 
if and only if both~$\omega$ and~$\Upsilon$ vanish 
when pulled back to~$L$.

\subsection{The group action}\label{ssec: son}
The group~$\SO(n)$ will be taken to act on~$\bbC^{n+1}
=\bbC\oplus\bbC^n$ in the manner as described in the 
introduction, namely, trivially on the first~$\bbC$-summand
and as the $\bbC$-linear extension of its standard action 
on~$\bbR^n$ in the remaining~$\bbC^n$.

\subsection{The reductions}\label{ssec: reductions}
The momentum mapping of this action~$\mu:\bbC^{n+1}\to\euso(n)$ 
is then (up to a constant scale factor that will be irrelevant
in what follows)
\begin{equation}
\mu(z_0,\ldots,z_n) = (x_iy_j-y_ix_j)_{1\le i,j\le n} \in\euso(n).
\end{equation}

\subsubsection{Nonzero momentum in the case $n=2$}

When~$n>2$, the algebra~$\euso(n)$ has trivial center, 
so, by the remarks in~\S\ref{sssec: genG}, there is only 
need to consider~$\mu^{-1}(0)$, which is the case of most 
concern in this article.

However, when~$n=2$, the algebra~$\euso(2)\simeq\bbR$ 
is abelian and the level sets of~$\mu$ are the
hypersurfaces~$H_c$ of the form~$x_1y_2-y_1x_2 = c$,
where~$c\in\bbR$ is any constant. 

When~$c\not=0$, the $\mu$-level set~$H_c\subset\bbC^3$ is 
smooth and~$\SO(2)\simeq S^1$ acts freely on~$M_c$.  
The reduced space~$M_c = \SO(2)\backslash H_c$ is thus
a smooth $4$-manifold endowed with a nonintegrable almost complex
structure~$J_c$ such that the $J_c$-complex curves in~$M_c$
are the $\SO(2)$-quotients of the~$\SO(2)$-invariant special 
Lagrangian $3$-folds in~$\bbC^3$ that lie in the level 
set~$H_c$.  

Each of the~$(M_c,J_c)$ with~$c\not=0$ is equivalent 
to~$(M_1,J_1)$, so there is really only one case of
nonzero momentum that needs to be treated.  
In any case, since there are no singular $\SO(2)$-orbits 
involved, these cases are reduced to the study of complex 
curves in almost complex $4$-manifolds and will not be 
discussed further here.

\subsubsection{The zero momentum almost complex structure}
\label{sssec: zeromuAC}

Henceforth, to avoid having to continually mention trivial
cases, it will be assumed that~$n>1$.

The locus~$\mu^{-1}(0)\subset\bbC^{n+1}$ is singular, 
with its singular locus 
consisting of the fixed points of the~$\SO(n)$-action,
i.e., the points of the form~$(z_0,0,\ldots,0)$.
It can be parametrized via the mapping
$\Phi:\bbC\times\bbC\times S^{n-1}\to\mu^{-1}(0)$ 
defined by
\begin{equation}
\Phi(w,\zeta,{\bf u})
= (w,\zeta\,u_1,\zeta\,u_2,\ldots,\zeta u_n)
\end{equation}
where~$w,\zeta\in\bbC$ and~${\bf u} = (u_1,\ldots,u_n)\in S^{n-1}$.
Note that~$\Phi(w,\zeta, {\bf u}) = \Phi(w,\, -\zeta,\,-{\bf u})$,
so that~$\Phi$ is a $2$-to-$1$ diffeomorphism away from the 
locus~$\zeta=0$, which is mapped by~$\Phi$
onto the fixed locus of the $\SO(n)$-action.

Computation yields
\begin{equation}
\Phi^*(\omega) 
= {\ts\frac\iC2}\bigl(\d w\w \d\bar w +\d\zeta\w\d\bar\zeta\bigr)
\end{equation}
and
\begin{equation}
\Phi^*(\Upsilon)  
= {\ts\frac1n}\text{Im}\bigl(\d w \w \d(\zeta^n)\bigr)\w \Omega
\end{equation}
where~$\Omega$ is the standard volume form on~$S^{n-1}$, i.e.,
\begin{equation}
\Omega = u_1\,\d u_2\w\ldots\w\d u_n +\cdots 
+(-1)^{n-1}u_n\,\d u_1\w\ldots\w\d u_{n-1}\,.
\end{equation}

It follows that any $\SO(n)$-invariant special 
Lagrangian~$L\subset\bbC^{n+1}$ (of momentum zero in case~$n=2$)
that does not meet the fixed locus is of the 
form~$L=\Phi\bigl(\Sigma\times S^{n-1}\bigr)$ for some
surface~$\Sigma\subset\bbC^2$ that does not meet the line~$\zeta=0$
and that is an integral manifold of the two $2$-forms
\begin{equation}
\begin{aligned}
\Psi_1 &= {\ts\frac\iC2}\bigl(\d w\w \d\bar w +\d\zeta\w\d\bar\zeta\bigr),\\
\Psi_2 &=  \text{Im}\bigl(\d w \w \d(\zeta^n)\bigr).
\end{aligned}
\end{equation}
Let~$M_0\subset\bbC^2$ denote the complement of the line~$\zeta=0$.
Then~$\Psi_1$ and~$\Psi_2$ are linearly independent on~$M_0$.

Moreover, on~$M_0$, the $2$-forms $\Psi_1$ and~$\Psi_2$ are 
multiples of the real and imaginary parts of the 
decomposable complex-valued $2$-form
\begin{equation}
\Psi = \omega_1\w\omega_2
= -\frac2{n|\zeta|^{n-1}}\,\Psi_2 + 2\iC\,\Psi_1
\end{equation}
where~$\omega_1$ and~$\omega_2$ are the $1$-forms on~$M_0$ 
defined by the formulae
\begin{equation}\label{eq: omega12def}
\omega_1 = \d w 
+ \iC\,\frac{{\overline\zeta}^{n-1}\,\d\overline\zeta}{|\zeta|^{n-1}},
\qquad
\omega_2 = \d{\overline{w}}
+ \iC\,\frac{{\zeta}^{n-1}\,\d\zeta}{|\zeta|^{n-1}}.
\end{equation}
Since~$\omega_1\w\omega_2\w\overline{\omega_1}\w\overline{\omega_2}\not=0$,
it follows that there exists a unique almost complex structure~$J_0$ 
on~$M_0$ for which $\omega_1$ and~$\omega_2$ furnish a basis 
of~$\Omega^{1,0}(M_0,J_0)$.

\begin{remark}[Nonintegrability of~$J_0$]\label{rem: nonint}
Although it will not be needed in the rest of this article, the
reader might like to know a bit more about the almost complex
structure~$J_0$, so some information about it will be mentioned here. 

First, as is easily computed,~$J_0$ is not integrable.  In fact,
its Nijenhuis tensor is nowhere vanishing.

Second, the almost complex manifold~$(M_0,J_0)$ is homogeneous:  
The diffeomorphisms $F_{a,b}:M_0\to M_0$ defined by
\begin{equation}\label{eq: Fab}
F_{a,b}(w,\zeta) 
= \left(\,\frac{{\overline{b}}^n}{|b|^{n-1}}\,w + a,\,\,b\,\zeta\,\right)
\end{equation}
for constants~$a\in\bbC$ and~$b\in\bbC^*$ preserve 
the almost complex structure~$J_0$ and it is clear that these
mappings generate a group that acts simply transitively on~$M_0$.

It is not difficult to show that the mappings~$F_{a,b}$ as
defined in~\eqref{eq: Fab} together with the 
involution~$A:M_0\to M_0$
defined by
\begin{equation}\label{eq: A}
A(w,\zeta) 
= (\,\overline{w},\,\overline{\zeta}\,) 
\end{equation}
generate the group of automorphisms of~$(M_0,J_0)$.  In fact,
a slightly stronger statement is true:  Any local automorphism
of~$(M_0,J_0)$ defined on a connected open subset of~$M_0$
extends uniquely to a global automorphism and is either
of the form~$F_{a,b}$ or~$A{\circ}F_{a,b}$.  

This last claim follows from the properties of the $(1,0)$-forms
\begin{equation}\label{eq: eta12def}
\eta_1 = \frac{|\zeta|^{n-1}}{{\overline\zeta}^n}\,\d w 
+ \iC\,\frac{\d\overline\zeta}{\overline\zeta}\,,
\qquad
\eta_2 = \frac{|\zeta|^{n-1}}{{\zeta}^n}\,\d{\overline{w}}
+ \iC\,\frac{\d\zeta}{\zeta}\,,
\end{equation}
which are invariant under the action of~$F_{a,b}$ and
are exchanged by~$A$.  

Inspection shows that the mapping~$C:M_0\to M_0$ defined by
\begin{equation}\label{eq: Cdef}
C(w,\zeta)=(w,\E^{\iC\pi/n}\zeta)
\end{equation} 
is $J_0$-antilinear (since~$C^*\omega_1 = \overline{\omega_2}$ 
and~$C^*\omega_2 = \overline{\omega_1}$),
a fact that will be useful below.

Third, the almost complex structure~$J_0$ on~$M_0$ cannot 
be extended continuously across the line~$\zeta=0$.  
However, in view of the fact that this line is the fixed
locus of the `conjugation'~$C$, one can think of the
line~$\zeta=0$ as a sort of `singular' totally real
submanifold of~$(\bbC^2,J_0)$.
\end{remark}

\begin{remark}[The double cover]\label{rem: dblcvr}
The reader may have noticed that~$M_0$ 
is not the symplectic quotient~$\SO(n)\backslash\mu^{-1}(0)^*$ 
but rather is a double cover of it.  In fact, when~$n$ is odd,
it is necessary to take this double cover
since~$\mu^{-1}(0)^*$ is not orientable when~$n$ is odd. 
Consequently, in this case,~$J_0$ is only defined 
up to a sign on~$\SO(n)\backslash\mu^{-1}(0)^*$.

There are other reasons for working on~$M_0$.  As will be
seen below, a nonsingular $\SO(n)$-invariant special Lagrangian
submanifold~$L\subset\bbC^{n+1}$ that meets the fixed locus
will be represented by a smooth surface~$\Sigma\subset M_0$
that extends across the line~$\zeta=0$ to a smoothly embedded 
surface in~$\bbC^2$ that meets the line~$\zeta=0$ in a
smooth analytic arc.
\end{remark}
 
\subsection{Invariant special Lagrangian planes}
\label{ssec: invsLagplanes}
Suppose that a nonsingular $\SO(n)$-invariant special 
Lagrangian~$L\subset\bbC^{n+1}$ 
meets the fixed locus~$\bbC$ of the $\SO(n)$-action 
at a point~$z\in L$.  Then the tangent plane~$T_zL$
must be a $\SO(n)$-invariant special Lagrangian $(n{+}1)$-plane.

There is a circle of such~$(n{+}1)$-planes:  
For each~$\psi$, there is the $(n{+}1)$-plane~$P_\psi$ 
defined by the linearly independent equations
\begin{equation}\label{eq: Ppsidef}
\begin{aligned}
0 &=  \cos n\psi\,\d y_0 + \sin n\psi\,\d x_0\\
  &=  \cos\psi\,\d y_1 - \sin\psi\,\d x_1  = \cdots
   = \cos\psi\,\d y_n - \sin\psi\,\d x_n\,.
\end{aligned} 
\end{equation}
Of course~$P_\psi=P_{\psi+\pi}$, but the planes~$\{P_\psi\ \vrule\ 
0\le\psi<\pi\}$ are pairwise distinct.  

Note that each~$P_\psi$ intersects~$\bbC$ in a real $1$-dimensional
linear subspace; that, conversely, each (real) $1$-dimensional 
linear subspace of~$\bbC$ lies in exactly $n$ of these $\SO(n)$-invariant 
special Lagrangian $(n{+}1)$-planes; and that the projections of
these $n$ special Lagrangian $(n{+}1)$-planes into~$\bbC^n$ 
are pairwise disjoint.

\subsubsection{The $\lambda$-action}\label{sssec: lambdaact}
Let~$\lambda = \E^{\pi\iC/n}$, so that~$\lambda$ generates
a multiplicative cyclic subgroup of order~$2n$,
denoted~$\bbZ_{2n}\subset S^1\subset\bbC$.

The $\bbZ_{2n}$-action on~$\bbC^{n+1}$ defined by
\begin{equation}
\lambda^j{\star}(z_0,,z_1,\ldots,z_n) 
= \bigl(z_0,\,\lambda^j z_1,\,\ldots,\,\lambda^j z_n\bigr)
\end{equation}
pulls back the holomorphic volume form~$\d z$ to~$(-1)^j\d z$
and commutes with the~$\SO(n)$-action.
It follows that this action 
carries special Lagrangian $(n{+}1)$-planes
to other special Lagrangian $(n{+}1)$-planes
(but may reverse their orientations)
and permutes the~$n$ $\SO(n)$-invariant special Lagrangian 
$(n{+}1)$-planes that contain a given fixed (real) 
line in~$\bbC$.  

In particular, if~$L\subset\bbC^{n+1}$ is special Lagrangian, 
then $\lambda^j{\star}L$ is also special Lagrangian
for~$0\le j< 2n$ and satisfies~$\lambda^j{\star}L\cap\bbC 
= L\cap\bbC$.  Note that~$\lambda^{j+n}{\star}L$ 
and~$\lambda^{j}{\star}L$ are tangent along their common
intersection with~$\bbC$.  As will be seen below
in Remark~\ref{rem: ln=id}, these two special Lagrangian 
submanifolds are actually equal in a neighborhood of~$\bbC$.

\subsubsection{A commuting action}\label{sssec: comact} 
The~$\SO(n)$-action commutes with the group of
special Lagrangian symmetries of~$\bbC^{n+1}$ of the form
\begin{equation}\label{eq: circ_act}
\Phi_{a,\theta}(z_0,z_1,\ldots,z_n) 
= \bigl(\E^{n\iC\theta}\,z_0+a,\,
\E^{-\iC\theta}\,z_1,\,\ldots,\,\E^{-\iC\theta}\,z_n,
         \bigr).
\end{equation}
where~$a\in\bbC$ and~$\theta\in\bbR/(2\pi\bbZ)$ are constants.
Note that this group action (which \emph{does} 
preserve the special Lagrangian calibration)
acts transitively on the $(n{+}1)$-planes of the form~$P_\psi$.  
Also, this group action preserves the fixed locus~$\bbC$ and 
acts on~$\bbC$ as the group of Euclidean isometries of~$\bbC$
(regarded as a real $2$-plane). 

\subsection{The fixed locus}\label{ssec: fixL}

The goal of this subsection is to examine the geometry
of an $\SO(n)$-invariant special Lagrangian submanifold~$L\subset\bbC^n$
that meets the fixed locus~$\bbC\subset\bbC^{n+1}$
then apply the results to prove the following 
(which was proved by Dominic Joyce in the case~$n=2$):

\begin{proposition}\label{prop: Anonsing}
Suppose that~$L\subset\bbC^{n+1}$ is an embedded nonsingular
special Lagrangian submanifold that is $\SO(n)$-invariant
and that meets the fixed locus~$\bbC$.  Then~$A = L\cap\bbC$
is an embedded nonsingular real-analytic curve in~$\bbC$.
If~$n>2$, then~$A$ has no compact component.
\end{proposition}

\begin{proof}
Let~$L\subset\bbC^{n+1}$ be an embedded, nonsingular 
$\SO(n)$-invariant special Lagrangian submanifold 
that meets the fixed line~$\bbC\subset\bbC^{n+1}$ 
at a point~$z\in L\cap\bbC$. 

After applying an action of the form~\eqref{eq: circ_act}, 
it can be assumed that~$z=0$ and that~$T_0L = P_0 = \bbR^{n+1}$. 

Since~$L$ is embedded and Lagrangian 
and since~$T_0L=\bbR^{n+1}$, 
it follows that, in some neighborhood of~$0\in\bbC^{n+1}$, 
the submanifold~$L$ can be parametrized in the form
\begin{equation}\label{eq: F_graph}
\left(
x_0+\iC\,\frac{\partial F}{\partial x_0}(x_0,x_1,\ldots,x_n),\,
\ldots,\,
x_n+\iC\,\frac{\partial F}{\partial x_n}(x_0,x_1,\ldots,x_n),\,\right)
\end{equation}
for some function~$F$ that is defined on a neighborhood 
of~$0\in\bbR^{n+1}$ and has all of its first and second partials 
vanishing there.  The function~$F$ can be made unique by requiring 
that~$F(0,\ldots,0)=0$, so assume this.  

Since~$L$ is nonsingular and special Lagrangian (and hence minimal), 
the known regularity of minimal submanifolds~\cite{MR21:5070} 
implies that~$L$ is real-analytic and hence that~$F$ 
is real-analytic also.

Because~$F$ is invariant under the action of~$\SO(n)$, 
there exists a real-analytic function~$\phi$ defined 
in a neighborhood~$V\subset\bbR^2$ 
of~$(0,0)$ that is even in the second variable (i.e.,
$\phi(t,\sigma)=\phi(t,-\sigma)$) and that satisfies
\begin{equation}\label{eq: F_as_phi}
F(x_0,x_1,\ldots,x_n) 
= \phi\left(x_0,\,\sqrt{{x_1}^2{+}\cdots{+}{x_n}^2}\,\right). 
\end{equation}
for~$(x_0,\ldots,x_n)$ sufficiently near the origin in~$\bbR^{n+1}$.
(The reason for defining~$\phi$ as an even function of~$\sigma=
\sqrt{{x_1}^2{+}\cdots{+}{x_n}^2}$ rather than directly as a
function of~${x_1}^2{+}\cdots{+}{x_n}^2$ is that it leads to
a more manageable equation in the uniqueness analysis to be
done below.)  

Because~$\phi$ is even in its second argument, 
the quotient~$\phi_\sigma(t,\sigma)/\sigma$ 
is a real-analytic function (also even in~$\sigma$) 
on~$V\subset\bbR$.  Thus, the graph of~$F$ 
can be parametrized analytically near~$0\in\bbC^{n+1}$ 
in the form
\begin{equation}\label{eq: phi_graph}
\left(
t+\iC\,\phi_t(t,|x|),\,
x_1+\iC\,x_1\frac{\phi_\sigma(t,|x|)}{|x|},\,
\ldots,\,
x_n+\iC\,x_n\frac{\phi_\sigma(t,|x|)}{|x|}\,\right)
\end{equation}
where~$|x|=\sqrt{{x_1}^2{+}\cdots{+}{x_n}^2}$.
The condition that~$T_0L=\bbR^{n+1}$ is then 
equivalent to the conditions that~$\phi_t(t,\sigma)$
have vanishing differential at~$(t,\sigma)=(0,0)$
and that the smooth function~$\phi_\sigma(t,\sigma)/\sigma$ 
vanish at~$(t,\sigma)=(0,0)$.

An immediate consequence of the representation~\eqref{eq: phi_graph} 
is that, in a neighborhood of the origin,~$L\cap\bbC$ 
consists of the points of the form
\begin{equation}
\left(\,t+\iC\,\phi_t(t,0),\,0,\ldots,0\right)
\end{equation}
for~$|t|$ sufficiently small, which is a nonsingular
real-analytic curve. 

Since~$z$ was an arbitrary point of~$L\cap\bbC$, it follows 
that $L\cap\bbC$ is a nonsingular embedded real-analytic 
curve~$A\subset \bbC$.

Finally, suppose that~$A= L\cap\bbC$ has a compact component,
i.e., an embedded closed curve~$A_0\subset A$.  Choose a
periodic, nonsingular parametrization~$(x,y):\bbR\to A_0$
with period~$1$, i.e., $x(t{+}1)=x(t)$ and~$y(t{+}1)=y(t)$.
There will then exist a smooth function~$\theta:\bbR\to\bbR$
such that
\begin{equation}
x'(t) = \cos\theta(t)\,\sqrt{x'(t)^2+y'(t)^2},
\qquad
y'(t) = -\sin\theta(t)\,\sqrt{x'(t)^2+y'(t)^2}.
\end{equation}
Since~$A_0$ is an embedded curve, it has rotation number~$\pm1$,
i.e.,~$\theta(t{+}1) = \theta(t)\pm2\pi$.  By reversing the
orientation of the parametrization if necessary, it can be 
supposed that~$\theta(t{+}1) = \theta(t)+2\pi$.

Now, the special Lagrangian planes that contain the tangent line
to~$A$ at~$z = \bigl(x(0),y(0)\bigr)$ are, by construction, of
the form~$P_{\theta(0)/n+k\pi/n}$ for~$k=0,1,\ldots,n{-}1$, so
one of these is~$T_zL$.  Fix~$k$ so that~$T_zL=P_{\theta(0)/n+k\pi/n}$.
Then, by continuity, it follows that
\begin{equation}
T_{\bigl(x(t),y(t)\bigr)}L = P_{\theta(t)/n+k\pi/n}
\end{equation}
for all~$t$.  However, since~$x$ and~$y$ are periodic of
period~$1$, it would then follow that
\begin{equation}
\begin{aligned}
P_{\theta(0)/n+k\pi/n} &= T_{\bigl(x(0),y(0)\bigr)}L
= T_{\bigl(x(1),y(1)\bigr)}L\\
&= P_{\theta(1)/n+k\pi/n}= P_{\theta(0)/n+2\pi/n+k\pi/n}.
\end{aligned}
\end{equation}
However, if~$n$ were greater than~$2$, the first and last 
of the planes in this string of equalities could not
be equal.  Thus, $n=2$, as claimed.
\end{proof}

\begin{example}[Extending the unit circle]\label{ex: extunitcirc}
The last statement of Proposition~\ref{prop: Anonsing}
may seem surprising at first.  However, consider the
case in which~$A$ is the unit circle in~$\bbC\subset\bbC^{n+1}$.
Using the circle invariance to reduce the integration 
problem to an \textsc{ode} (as in Example~\ref{ex: cohom1}), 
one finds that the locus
\begin{equation}
L = \left\{\, (z,\zeta {\bf u})\ \vrule\ |\zeta|^2 = n\bigl(|z|^2-1\bigr),
\ \text{Re}(z\zeta^n) = 0,\,{\bf u}\in S^{n-1}\,\right\}
\subset\bbC^{n+1}
\end{equation}
is smooth away from~$A = L\cap\bbC$ and is special
Lagrangian.  Near each point of~$A$, the locus~$L$ is
a union of $n$ smooth distinct sheets intersecting 
pairwise in~$A$.  

When~$n$ is odd,~$L\setminus A$ is connected, so~$L$ 
is analytically irreducible.
When~$n$ is even,~$L\setminus A$ has two connected components, 
and, in fact, $L$ is the union of two distinct
analytically irreducible pieces: $L_+$, 
on which~$\text{Im}(z\zeta^n)$ is nonnegative, 
and~$L_-$, on which~$\text{Im}(z\zeta^n)$ is nonpositive.
Near each point of~$A$, each of~$L_+$ and~$L_-$ is
a union of~$\frac12n$ smooth distinct sheets.
Only when~$n=2$ are~$L_+$ and~$L_-$ smooth.

This behavior for the unit circle is typical for embedded
closed curves in general, as will be seen.
\end{example}

\begin{remark}[Invariance under~$\lambda^n$]
\label{rem: ln=id}
The representation~\eqref{eq: phi_graph} also shows
that~$F$ is even in the variables~$x_1,\ldots,x_n$.  
In particular, this implies that, if~$L\subset\bbC^{n+1}$ 
is the graph of~$F$ regarded as a special Lagrangian 
submanifold, then~$\lambda^n{\star}L = L$ on some 
neighborhood of~$0\in\bbC^{n+1}$. 

By analytic continuation, it follows that, for any embedded, 
non-singular, $\SO(n)$-invariant, connected special 
Lagrangian~$L\subset\bbC^{n+1}$ that meets~$\bbC$,
the submanifolds~$\lambda^n{\star}L$ and~$L$ are
equal in some open neighborhood of~$L\cap\bbC$.
\end{remark}

\begin{remark}[Image in~$M_0$]\label{rem: im_inM0}
It is a consequence of the proof that 
there is a neighborhood~$U\subset\bbC^{n+1}$
of~$0$ such that~$U\cap L = \Phi(\Sigma_\phi\times S^{n-1})$ where
$\Sigma_\phi\subset\bbC^2$ is the analytically embedded surface
\begin{equation}
\Sigma_\phi = \left\{\,\bigl(\,t+\iC\,\phi_t(t,\sigma),\,
                             \sigma+\iC\,\phi_\sigma(t,\sigma)\,\bigr)
    \ \vrule\ (t,\sigma)\in V\ \right\}.
\end{equation}

	Note that, when~$n$ is odd, 
the embedding~$\iota_\phi:V\to\bbC^2$ defined by
\begin{equation}
\iota_\phi(t,\sigma) =  \bigl(\,t+\iC\,\phi_t(t,\sigma),\,
                             \sigma+\iC\,\phi_\sigma(t,\sigma)\,\bigr)
\end{equation}
pulls back~$\omega_1$ (which is only defined on~$M_0$ via
the formulae~\eqref{eq: omega12def}) to a complex-valued $1$-form 
that extends smoothly across the curve~$\sigma=0$.  This is
not so surprising since, when~$n$ is odd, the mapping~$(w,\zeta)
\mapsto(w,-\zeta)$ is $J_0$-antilinear on~$M_0$ and $\iota_\phi$  
intertwines this mapping with the orientation reversing mapping
$(t,\sigma)\mapsto(t,-\sigma)$ on~$V$.  In any case, 
~${\iota_\phi}^*(\omega_1)$ is a $(1,0)$-form for a natural 
complex structure on~$V$ that makes~$\iota_\phi$ into a 
$J_0$-complex curve away from the locus~$\sigma=0$.  

The picture when~$n$ is even is slightly more complicated
and is left to the reader.
\end{remark}

\section{Local uniqueness}\label{sec: locallyunique}

Suppose now that~$A\subset \bbC$ is 
a connected, nonsingular real-analytic curve.
The goal now is to determine whether~$A$ is the fixed locus 
of a nonsingular $\SO(n)$-invariant special Lagrangian 
$(n{+}1)$-fold~$L\subset\bbC^{n+1}$ and, if so, in
how many ways.

\begin{remark}[Lack of uniqueness]\label{rem: nonunique}
As has already been remarked, if~$A = L\cap\bbC$
for some embedded nonsingular $\SO(n)$-invariant 
special Lagrangian $(n{+}1)$-fold~$L$, 
then~$A=\lambda^j{\star}L\cap\bbC$ 
for any integer~$j$ in the range~$0\le j < n$.
Moreover, since~$L$ is embedded, $A$ has an open 
neighborhood $U\subset\bbC^{n+1}$ so that
\begin{equation}
\lambda^j{\star}L \cap \lambda^k{\star}L \cap U = A
\end{equation}
for any~$j$ and~$k$ satisfying~$0\le j < k < n$.
 
One consequence of the analysis to be done
below is that any embedded nonsingular
special Lagrangian $(n{+}1)$-fold~$L'\subset\bbC^{n+1}$ 
that contains~$A$ agrees with~$\lambda^j{\star}L$ in
some open neighborhood of~$A$ in~$\bbC^{n+1}$ 
for some integer~$j$ in the range~$0\le j < n$.%
\footnote{Recall from Remark~\ref{rem: ln=id} 
that $\lambda^{j+n}{\star}L$ and $\lambda^j{\star}L$ 
agree in some open neighborhood of~$A$. Thus, one can 
restrict the range of~$j$ to~$0\le j<n$ as claimed.} 
This result will follow from general considerations once it
is shown that a local version of this uniqueness holds.
\end{remark}

\subsection{Reduction to an equation}\label{ssec: red2pde}
Using the action~\eqref{eq: circ_act}, to understand
the local picture, it suffices to understand
the case where the curve~$A$ passes through the origin~$0\in\bbC^{n+1}$, 
its tangent there is spanned by~$\p/\p x_0$, and~$T_0L=\bbR^{n+1}$ is
spanned by~$\p/\p x_0,\,\p/\p x_1,\,\ldots,\,\p/\p x_n$.

As in~\S\ref{ssec: fixL}, it follows that, in a neighborhood 
of~$0\in\bbC^{n+1}$, $L$ can be described as a graph 
of the form~\eqref{eq: F_graph} 
for some function~$F$ of the form~\eqref{eq: F_as_phi}, 
where~$\phi$ is a real-analytic function on a neighborhood 
of~$0\in\bbR^2$.

By hypothesis,~$A$ can be parametrized 
near~$0\in\bbC^{n+1}$ in the form
\begin{equation}\label{eq: A_graph}
\bigl(t+\iC\,f'_0(t),0,\ldots,0\bigr)
\end{equation}
for some function~$f_0$ that is real-analytic 
in a neighborhood of~$0\in\bbR$
and abides~$f'_0(0)= f''_0(0)=0$.  
The function~$f_0$ can be made unique
by requiring that~$f_0(0)=0$.  
(The reason for starting with~$f'_0$ 
instead of~$f_0$ should be apparent.)

Now, the condition that~$L$ contain~$A$ becomes 
the condition
\begin{equation}\label{eq: phi_init0}
\phi(t,0)=f_0(t).
\end{equation}
(Bear in mind the normalizations~$F(0,\ldots,0)=f_0(0)=0$.)
Furthermore, the condition that~$T_0L = \bbR^{n+1}$ implies 
\begin{equation}\label{eq: phi_init2}
\phi_{\sigma\sigma}(0,0) = 0.
\end{equation}

Let~$\sigma$ stand for~$\sqrt{{x_1}^2+\cdots+{x_n}^2}$ 
and let~$t$ stand for~$x_0$
as the coordinates in the domain of~$\phi$ in~$\bbR^2$.  
The condition that the $1$-graph
of~$F$ as defined in~\eqref{eq: F_graph} be special Lagrangian is, 
of course, a second order partial differential equation on~$F$.  
This equation can be expressed 
in terms of~$\phi$ in the form
\begin{equation}
\text{Im}\left((\sigma+\iC\,\phi_\sigma)^{n-1}\,
\d(\sigma+\iC\,\phi_\sigma)\w\d(t+\iC\,\phi_t)\right) = 0,
\end{equation}
or, in more classical \textsc{pde} terms:
\begin{equation}\label{eq: phi_pde}
\text{Im}\left((\sigma+\iC\,\phi_\sigma)^{n-1}
\bigl((1+\iC\,\phi_{tt})(1+\iC\,\phi_{\sigma\sigma})
+{\phi_{\sigma t}}^2\bigr)\right) = 0,
\end{equation}
which is a singular, second order Monge-Amp\`ere equation 
that is elliptic when~$\sigma\not=0$ but degenerate 
along~$\sigma=0$.

\begin{proposition}\label{prop: phiunique}
Let~$f_0$ be a real-analytic function defined on an open interval 
containing~$0\in\bbR$ and that satisfies~$f_0(0)=f'_0(0)=f''_0(0)=0$.
Then there is at most one real-analytic function~$\phi$ defined
on a neighborhood of~$(0,0)\in\bbR^2$ that satisfies~$\phi(t,\sigma)
=\phi(t,-\sigma)$, the equation~\eqref{eq: phi_pde}, and 
the initial conditions~\eqref{eq: phi_init0} and~\eqref{eq: phi_init2}.
\end{proposition}

\begin{proof}
Any such~$\phi$ must have a power series expansion of the form
\begin{equation}\label{eq: phi_powsers}
\phi(t,\sigma) = f_0(t) 
+ {\ts\frac1{2!}}f_1(t)\sigma^2 
+ {\ts\frac1{4!}}f_2(t)\sigma^4 
+ {\ts\frac1{6!}}f_3(t)\sigma^6 + \cdots\ .
\end{equation}
where the~$f_i$ for~$i>0$ are real-analytic on some interval~$|t|\le\tau$ 
and satisfy a bound of the form~$|f_i(t)|\le C/M^i$
when~$|t|<\tau$ for some constants~$C>0$ and~$M>0$.  Moreover,
by~\eqref{eq: phi_init2}, $f_1$ must satisfy~$f_1(0)=0$.

To prove uniqueness, it suffices to show that the
equation~\eqref{eq: phi_pde} together with the specified 
initial conditions determine~$f_i$ uniquely for~$i>0$. 
Write~\eqref{eq: phi_pde} in the equivalent form
\begin{equation}\label{eq: phi_pde2}
\text{Im}\left(\Bigl(1+\iC\,\frac{\phi_\sigma}{\sigma}\Bigr)^{n-1}
\bigl((1+\iC\,\phi_{tt})(1+\iC\,\phi_{\sigma\sigma})
+{\phi_{\sigma t}}^2\bigr)\right) = 0,
\end{equation}

Now, substituting~\eqref{eq: phi_powsers} 
into~\eqref{eq: phi_pde2}, collecting like powers of~$\sigma$,
and considering the coefficient of~$\sigma^0$ yields
\begin{equation}\label{eq: f1constraint}
\text{Im}\left(\bigl(1+\iC\,f_1(t)\bigr)^n
               \bigl(1+\iC\,f''_0(t)\bigr)\right) = 0.
\end{equation}
This is an algebraic equation for~$f_1$ in terms of~$f''_0$
that has up to~$n$ distinct roots.  However, by hypothesis
and initial condition,~$f''_0(0)=f_1(0)=0$, so there is only 
one continuous choice of~$f_1$ that will satisfy the given
initial condition, namely
\begin{equation}\label{eq: f1formula}
f_1(t) = -\tan\left(\frac{\tan^{-1}\bigl(f''_0(t)\bigr)}{n}\right).
\end{equation}
Thus, assume henceforth that~$f_1$ is defined 
by~\eqref{eq: f1formula}.  For simplicity, set
\begin{equation}
R(t)  = \bigl(1+\iC\,f_1(t)\bigr)^n
               \bigl(1+\iC\,f''_0(t)\bigr)
\end{equation}
and note that~$R$ is real-valued and satisfies~$R(0)=1$.

Let~$\tau>0$ be such that~$f_0$, $f_1$, and~$1/R$ have
convergent power series in the interval~$-\tau<t<\tau$.   

Now, the derivatives of~$\phi$ that appear
on the left hand side of~\eqref{eq: phi_pde2} have 
convergent power series expansions in~$\sigma$ of the form
\begin{equation}\label{eq: phidersers}
\begin{aligned}
\phi_{tt}(t,\sigma) 
&=  \sum_{k=0}^\infty{\frac1{(2k)!}}f''_k(t)\sigma^{2k},
&\phi_{\sigma t}(t,\sigma) 
&=  \sum_{k=0}^\infty{\frac1{(2k{+}1)!}}f'_{k+1}(t)\sigma^{2k+1},\\
\phi_{\sigma\sigma}(t,\sigma) 
&=  \sum_{k=0}^\infty{\frac1{(2k)!}}f_{k+1}(t)\sigma^{2k},
&\frac{\phi_\sigma(t,\sigma)}{\sigma}
&=  \sum_{k=0}^\infty{\frac1{(2k{+}1)!}}f_{k+1}(t)\sigma^{2k}.\\
\end{aligned}
\end{equation}
Using these expressions, the formula for~$f_1$ and~$R$, and the definitions 
given above, it follows that the left hand side of~\eqref{eq: phi_pde2}
has a series expansion in~$\sigma$ of the form
\begin{equation}\label{eq: phi_pde2LHSsers}
\sum_{k=1}^\infty
\left[
\frac{R}{1+{f_1}^2}\,\frac{2k{+}n}{(2k{+}1)!}\,f_{k+1}
-Q_k(f_1,\ldots,f_k,f'_1,\ldots,f'_k,f''_0,\ldots,f''_k)
\right]\sigma^{2k}
\end{equation}
where each~$Q_k$ is an explicit polynomial of total degree
at most~$n{+}1$ in its~$3k{+}1$ arguments as listed.
(The essential points to note are: First, the coefficient 
of~$\sigma^0$ has already been set to zero by the definition 
of~$f_1$. Second, as the terms on the right hand sides 
of~\eqref{eq: phidersers} show, for~$k>0$, the coefficient 
of~$\sigma^{2k}$ in the $\sigma$-expansion
of the left hand side of~\eqref{eq: phi_pde2} is a sum of 
products of terms whose coefficients only involve the quantities
$f_1,\ldots,f_k,f_{k+1},f'_1,\ldots,f'_k,f''_0,\ldots,f''_k$.
Moreover, since terms involving~$f_{k+1}$ in this sum 
can only come from terms that have a factor of~$f_{k+1}\sigma^{2k}$, 
the remaining factors in such a term must occur
as coefficients of~$\sigma^0$.  These are easily collected
and computed, leading to the expression given 
in~\eqref{eq: phi_pde2LHSsers}.)

Consequently, since~$f_0$ is given and~$f_1$ is
determined by~\eqref{eq: f1formula}, the functions~$f_{k+1}$ 
for~$k\ge1$ are determined recursively by the equations
\begin{equation}\label{eq: frecursion}
f_{k+1}
= \frac{1+{f_1}^2}{R}\,\frac{(2k{+}1)!}{2k{+}n}\,
Q_k(f_1,\ldots,f_k,f'_1,\ldots,f'_k,f''_0,\ldots,f''_k).
\end{equation}
In particular, each~$f_k$ has a convergent power series
on the interval~$-\tau<t<\tau$.

Thus, there is a unique formal power series 
solution~$\phi(t,\sigma)$ to~\eqref{eq: phi_pde} 
that is even in~$\sigma$ and satisfies the initial 
conditions~\eqref{eq: phi_init0} and \eqref{eq: phi_init2}, 
as was to be shown.
\end{proof}

\section{Local existence}\label{sec: localexist}

While the uniqueness theorem above demonstrates the
existence of a formal power series solution 
to~\eqref{eq: phi_pde} that satisfies the appropriate
initial conditions, proving that the formal series
converges on some neighborhood of~$(t,\sigma)=(0,0)$ 
in~$\bbR^2$ is not easy to do directly.  The crude
argument used to derive 
the recursion formula~\eqref{eq: frecursion} does
not give sufficient detail about the polynomials~$Q_k$
to allow any effective estimates to be done on the 
growth of the terms~$f_{k}$ as~$k$ tends to~$\infty$.

Moreover, because of the singular nature of the equations
involved, a direct appeal to the Cauchy-Kowalewski theorem
does not seem to be feasible.

\subsection{An existence theorem of G\'erard and Tahara}
However, by making use of a more subtle application of
the method of majorants, G\'erard and Tahara have proven
an existence and uniqueness theorem in the holomorphic
category that suffices to prove that the above series
(which is the only formal power series solution) does,
in fact, converge.  Their main existence result
can be found as Theorem~8.0.3 of their book~\cite{MR2001c:35056}
and concerns the existence of holomorphic solutions of
holomorphic singular partial differential equations.  
For the convenience of the reader, this result will
now be summarized in the case of second order equations,
which is all that will be needed for this article.

Their existence theorem applies to certain singular partial
differential equations for a function~$u = u(t,\sigma)$ 
of the form
\begin{equation}\label{eq: Geq}
G\bigr(\,t,\,\sigma,\,u,\,u_t,\,\sigma u_\sigma,\,u_{tt}, 
\,\sigma u_{\sigma t},\,\sigma^2 u_{\sigma\sigma}\,\bigl)=0,
\end{equation}
where~$G$ is a real-analytic function of its eight arguments
in a neighborhood of the origin in~$\bbR^8$.  

To avoid confusion in the statement of their result, it will
be important to adopt a notation that clearly distinguishes
the arguments and suggests their meanings.  For this reason,
the eight arguments of~$G$ will be written as
\begin{equation}\label{eq: Gvars}
G\bigl(t,\sigma, Z_{0,0}, Z_{0,1}, Z_{1,0}, Z_{0,2},
 Z_{1,1}, Z_{2,0}\bigr).
\end{equation}
The idea is that~$Z_{k,l}$ is the name of the variable
into which the derivative expression
\begin{equation}
\sigma^k\frac{\p^{k+l}u}{\p^k\sigma\,\p^lt}
\end{equation}
is to be substituted in~\eqref{eq: Gvars} 
in order to form the left hand side of~\eqref{eq: Geq}  
Also, in order to save writing, for any analytic
function on a neighborhood of the origin in~$\bbR^8$,
say
\begin{equation}
H\bigl(t,\sigma, Z_{0,0}, Z_{0,1}, Z_{1,0}, Z_{0,2},
 Z_{1,1}, Z_{2,0}\bigr),
\end{equation}
the expression~$H(t,{\bf0})$
will be an abbreviation for~$H(t,0,0,0,0,0,0,0)$.

\begin{theorem}[cf. Theorem~8.0.3 of~\cite{MR2001c:35056}]
\label{thm: G&T}
Let~$G=G\bigl(t,\sigma, Z_{0,0}, Z_{0,1}, Z_{1,0}, Z_{0,2},
 Z_{1,1}, Z_{2,0}\bigr)$ be a real-analytic function
defined on a neighborhood of the origin in~$\bbR^8$.
Suppose that~$G$ satisfies the following conditions:
\begin{enumerate}
\item $G(t,{\bf0}) \equiv 0$,
\vspace{3pt}
\item ${\ds\frac{\p G\hfil}{\p Z_{0,1}}}(t,{\bf0})\equiv
       {\ds\frac{\p G\hfil}{\p Z_{1,1}}}(t,{\bf0})\equiv
       {\ds\frac{\p G\hfil}{\p Z_{0,2}}}(t,{\bf0})\equiv0$,
\vspace{3pt}
\item ${\ds\frac{\p G\hfil}{\p Z_{2,0}}}(0,{\bf0})\not=0$,
\vspace{3pt}
\item $ {\ds\frac{\p G\hfil}{\p Z_{2,0}}}(0,{\bf0})\,k^2
       +{\ds\frac{\p G\hfil}{\p Z_{1,0}}}(0,{\bf0})\,k
       +{\ds\frac{\p G\hfil}{\p Z_{0,0}}}(0,{\bf0}) \not=0$
       for any integer~$k>0$.
\end{enumerate}
Then there is a unique real-analytic function~$u$
defined on an open neighborhood of~$(t,\sigma)=(0,0)\in\bbR^2$
that satisfies the equation
\begin{equation}
G\bigr(\,t,\,\sigma,\,u,\,u_t,\,\sigma u_\sigma,\,u_{tt}, 
\,\sigma u_{\sigma t},\,\sigma^2 u_{\sigma\sigma}\,\bigl)=0
\end{equation}
and the initial condition~$u(t,0)\equiv0$. 
\end{theorem}

\begin{remark}[The hypotheses of G\'erard and Tahara]
The reader may well wonder about the significance 
of the hypotheses listed above.  

The first two conditions are meant to make the initial
condition~$u(t,0)=0$ formally compatible with 
the equation~\eqref{eq: Geq}.  

The third condition ensures that the
equation~\eqref{eq: Geq} can be solved near the origin in~$\bbR^8$ 
for the expression~$\sigma^2u_{\sigma\sigma}$.  In fact, G\'erard 
and Tahara state their theorem for equations of the form
\begin{equation}
\sigma^2u_{\sigma\sigma} - 
F\bigr(\,t,\,\sigma,\,u,\,u_t,\,\sigma u_\sigma,
\,u_{tt},\,\sigma u_{\sigma t}\,\bigl)  =  0
\end{equation}
where~$F$ is a real-analytic function of its arguments 
on a neighborhood of the origin in~$\bbR^7$.  The form
in which it has been stated here is more convenient for
the intended application.

The fourth condition ensures that, when one tries to
solve~\eqref{eq: Geq} by substituting in a power series
of the form
\begin{equation}\label{eq: useries}
u(t,\sigma) = f_1(t)\,\sigma + f_2(t)\sigma^2 + \cdots  ,
\end{equation}
and setting equal to zero the coefficients of
the various powers of~$\sigma$ that result,
the resulting conditions on the~$f_i$ can be 
resolved into a recursion relation for the functions~$f_i$ 
without having to divide by any quantity that can vanish.
In other words, the fourth condition guarantees that there
will exist a formal power series solution of the form 
\eqref{eq: useries}.
\end{remark}

\begin{remark}[Extensions]
Theorem~8.0.3 of~\cite{MR2001c:35056} is considerably 
more general than the result stated here as 
Theorem~\ref{thm: G&T}.  In the first place, G\'erard and
Tahara consider equations
of arbitrary (finite) order, not just second order equations.
In the second place, their existence and uniqueness theorem
provides for the classification of more general solutions than
just the single-valued real-analytic solutions.  However, the
specialized version of their theorem stated above is all that 
will be needed in this article.
\end{remark}

\subsection{Application}
The results of G\'erard and Tahara can now be 
applied to show that the series solution~$\phi$ found in
the previous section converges.

\begin{proposition}\label{prop: phi_exist}
Let~$f_0$ be a real-analytic function defined on an open interval 
containing~$0\in\bbR$ and that satisfies~$f_0(0)=f'_0(0)=f''_0(0)=0$.
Then there exists a real-analytic function~$\phi$ defined
on a neighborhood of~$(0,0)\in\bbR^2$ that satisfies~$\phi(t,\sigma)
=\phi(t,-\sigma)$, the equation~\eqref{eq: phi_pde}, and 
the initial conditions~\eqref{eq: phi_init0} and~\eqref{eq: phi_init2}.
\end{proposition}  

\begin{proof}
It has already been shown that if the equation~\eqref{eq: phi_pde}
has a real-analytic solution~$\phi$ that is even in~$\sigma$
and satisfies the initial conditions~\eqref{eq: phi_init0} 
and \eqref{eq: phi_init2}, 
then its expansion in powers of
$\sigma$, must be given by~\eqref{eq: phi_powsers} where
$f_1$ is defined by~\eqref{eq: f1formula} and~$f_k$ for~$k\ge2$
are then determined recursively via~\eqref{eq: frecursion}.
Existence will follow once this series is shown to converge.

This suggests looking for a solution to~\eqref{eq: phi_pde} 
of the form
\begin{equation}\label{eq: phiasu}
\phi(t,\sigma) 
= f_0(t) + {\frac12}\bigl(f_1(t) + u(t,\sigma)\bigr)\sigma^2
\end{equation}
where~$f_1(t)$ is defined in terms of the real-analytic 
function~$f_0$ via~\eqref{eq: f1formula} and where~$u$ is
a real-analytic function defined on a neighborhood of the
origin in~$\bbR^2$ that satisfies~$u(t,0)\equiv0$. 

Now, substituting~\eqref{eq: phiasu} into the 
equation~\eqref{eq: phi_pde2} yields an
equation for~$u$ of the form~\eqref{eq: Geq} 
where~$G$ is taken to be the analytic function
\begin{equation}\label{eq: Gformula}
\begin{aligned}
G&=\text{Im}\Bigl\{
(1+\iC(f_1(t){+}2Z_{0,0}{+}Z_{1,0}))^{n-1} \Bigl[\sigma^2(f'_1(t){+}2Z_{0,1}{+}Z_{1,1})^2\\
&\qquad+\bigl(1+\iC(f_0''(t){+}{\ts\frac12}(f_1''(t){+}Z_{0,2})\sigma^2)\bigr)
\bigl(1{+}\iC(f_1(t){+}2Z_{0,0}{+}4Z_{1,0}{+}Z_{2,0})\bigr)\Bigr]\Bigr\}.
\end{aligned}
\end{equation}
In order to apply Theorem~\ref{thm: G&T}, the four hypotheses 
on~$G$ must now be verified.  

To verify the first hypothesis, note that
\begin{equation}
G(t,{\bf0}) 
= \text{Im}\bigl(\,(1+\iC f_1(t))^n(\,1 + \iC f_0''(t)\,) \bigr).
\end{equation}
Of course,~$f_1$ was defined so that the expression on the
right hand side would vanish identically.  

To verify the second hypothesis, 
note that each of the variables~$Z_{0,1}$,
$Z_{1,1}$ and~$Z_{0,2}$ occur in the formula~\eqref{eq: Gformula}
with a coefficient that is a positive power of~$\sigma$. 
Of course, this immediately implies that
\begin{equation}
{\ds\frac{\p G\hfil}{\p Z_{0,1}}}(t,{\bf0})\equiv
{\ds\frac{\p G\hfil}{\p Z_{1,1}}}(t,{\bf0})\equiv
{\ds\frac{\p G\hfil}{\p Z_{0,2}}}(t,{\bf0})\equiv0.
\end{equation}

The third and fourth hypotheses follow from the easily derived formulae
\begin{equation}
{\ds\frac{\p G\hfil}{\p Z_{2,0}}}(0,{\bf0}) = 1\,,\quad
{\ds\frac{\p G\hfil}{\p Z_{1,0}}}(0,{\bf0}) = n{+}3\,,\quad
{\ds\frac{\p G\hfil}{\p Z_{0,0}}}(0,{\bf0}) = 2n\,.
\end{equation}

Thus, Theorem~\ref{thm: G&T} applies:
There exists a unique real-analytic function~$u$ 
defined on a neighborhood of the origin in~$\bbR^2$ 
that satisfies the equation~\eqref{eq: Geq} where~$G$
is defined as in~\eqref{eq: Gformula} and the initial
condition~$u(t,0)\equiv0$.  

Because~$G$ as defined in~\eqref{eq: Gformula} is an even 
function of~$\sigma$, the function~$v$ defined by~$v(t,\sigma) 
= u(t,-\sigma)$ also satisfies~\eqref{eq: Geq}.  Since~$v(t,0)\equiv0$,
the uniqueness part of Theorem~\ref{thm: G&T} implies that
$v(t,\sigma)=u(t,\sigma)$, i.e., that~$u(t,-\sigma)=u(t,\sigma)$.

Finally, using this solution~$u$ to define~$\phi$ 
via~\eqref{eq: phiasu}, the desired local existence 
is established.  In particular, the series for~$\phi$
defined via~\eqref{eq: f1formula} and~\eqref{eq: frecursion} 
must be convergent on some neighborhood of~$(0,0)\in\bbR^2$.
\end{proof}

\section{Conclusions}\label{sec: conclusions}

In this last section, the local \textsc{pde} results
of the previous sections will be applied to prove
the main results.

The first result is an `$n$-uniqueness' theorem:

\begin{theorem}\label{thm: nuniqueL}
Suppose that~$L$ and ~$L'$ are nonsingular, 
embedded $\SO(n)$-invariant special Lagrangian
submanifolds of~$\bbC^{n+1}$ that have the
same fixed locus~$A\subset\bbC$ and that, 
moreover,~$A$ is connected.

Then there is an open $A$-neighborhood~$U\subset\bbC^{n+1}$
and a unique integer~$j$ satisfying~$0\le j < n$ such that
\begin{equation}
L'\cap U = \lambda^j{\star}L\cap U.
\end{equation}
Moreover, when~$0\le j < k < n$, the
submanifolds~$\lambda^j{\star}L\cap U$ and~$\lambda^k{\star}L\cap U$
intersect only along~$A$.
\end{theorem}

\begin{proof}
If~$A$ is empty, then there is nothing to prove,
so assume that~$A$ is nonempty 
and that~$z = (z_0,0,\ldots,0)\in A$.
Both~$T_zL$ and~$T_zL'$ are special Lagrangian
$(n{+}1)$-planes that contain the line~$T_zA$
and hence there is an integer~$j$ in the range~$0\le j < n$
and an angle~$\psi$ in the range~$0\le \psi < \pi$
such that~$T_zL = P_\psi$ and~$T_zL' = P_{\psi+(j/n)\pi}$.

By applying a motion in the group generated by the
mappings~\eqref{eq: circ_act}, it can be assumed that
$z=0\in\bbC^{n+1}$ and that~$T_zL=\bbR^{n+1}=P_0$.  
Then~$T_z(\lambda^{-j}{\star}L') = T_zL = P_0$ as well.  

Since~$L$ and~$L'$ are embedded, it follows that there 
is an open neighborhood~$W_z$ of~$0\in\bbR^{n+1}$ and
an open neighborhood~$U_z$ of~$0\in\bbC^{n+1}$ 
such that there exist unique real-analytic functions~$F$ 
and~$F'$ defined on~$W_z$ and vanishing at~$0\in\bbR^{n+1}$ 
so that
\begin{equation}
L\cap U_z = \left\{ 
\bigl(x_0+\iC\,{\ts\frac{\partial F}{\partial x_0}},\,
\ldots, x_n+\iC\,{\ts\frac{\partial F}{\partial x_n}}\bigr)
\ \vrule\ (x_0,\ldots,x_n)\in W_z\ \right\}
\end{equation}
while
\begin{equation}
\lambda^{-j}{\star}L'\cap U_z = \left\{ 
\bigl(x_0+\iC\,{\ts\frac{\partial F'}{\partial x_0}},\,
\ldots, x_n+\iC\,{\ts\frac{\partial F'}{\partial x_n}}\bigr)
\ \vrule\ (x_0,\ldots,x_n)\in W_z\ \right\}.
\end{equation}
Since~$T_z(\lambda^{-j}{\star}L') = T_zL = \bbR^{n+1}$,
it follows that~$F$ and~$F'$ have vanishing first and
second derivatives at~$0\in\bbR^{n+1}$.  As was argued
in the proof of Proposition~\ref{prop: Anonsing}, it 
follows that there unique exist real-analytic functions~$\phi$ 
and~$\phi'$ defined on a neighborhood of~$(0,0)$ in~$\bbR^2$ 
such that~$\phi(0,0)=\phi'(0,0)=0$ and such that
\begin{equation}
\begin{aligned}
F(x_0,x_1,\ldots,x_n) 
&= \phi\bigl(x_0,\sqrt{{x_1}^2{+}{\cdots}{+}{x_n}^2}\bigr)\\
F'(x_0,x_1,\ldots,x_n) 
&= \phi'\bigl(x_0,\sqrt{{x_1}^2{+}{\cdots}{+}{x_n}^2}\bigr)\\
\end{aligned}
\end{equation}
when~${x_0}^2{+}{\cdots}{+}{x_n}^2$ is sufficiently small.
Since~$L\cap U_z\cap \bbC = L'\cap U_z\cap \bbC = A \cap U_z$,
it follows that~$\phi(t,0)=\phi'(t,0)$ for~$|t|$ sufficiently
small.  Moreover, by construction,~$\phi$ and~$\phi'$ are
solutions of~\eqref{eq: phi_pde} that are even in the second 
variable (i.e., $\sigma$).  Because~$T_z(\lambda^{-j}{\star}L') 
= T_zL = \bbR^{n+1}$, it follows that~$\phi_{\sigma\sigma}(0,0)
=\phi'_{\sigma\sigma}(0,0)=0$ and, 
setting~$f_0(t)=\phi(t,0)=\phi'(t,0)$, that~$f_0(0)=f'_0(0)=f''_0(0)=0$.  

Now Proposition~\ref{prop: phiunique} implies that~$\phi=\phi'$,
which, in turn, 
implies that~$L\cap U_z = \lambda^{-j}{\star}L'\cap U_z$.

Since~$z$ was chosen arbitrarily in~$A$, it has been
established that every~$z\in A$ has an open neighborhood~$U_z
\subset\bbC^{n+1}$ such that there exists a unique 
integer~$j_z$ satisfying~$0\le j_z < n$ so 
that~$L'\cap U_z = (\lambda^{j_z}{\star}L) \cap U_z$.
Moreover, by shrinking~$U_z$ appropriately,
it can be arranged that, when~$0\le j < k \mod n$,
the two manifolds~$(\lambda^j{\star}L)$ 
and $(\lambda^k{\star}L)\cap U_z$ only intersect
along~$A\cap U_z$.  

In particular, it follows that the integer~$j_z$ is
locally constant in~$z$ and hence, since~$A$ is connected,
it follows that there is a unique integer~$j$
in the range~$0\le j<n$ such that~$j_z = j$ for all~$z\in A$.

Now let~$U$ be the union of the~$U_z$ as~$z$ ranges over~$A$.
Then~$U$ and~$j$ have the required properties.
\end{proof}

The second main result is an existence result:

\begin{theorem}\label{thm: existL}
Suppose that~$A\subset\bbC$ is an embedded, connected
real-analytic curve and that either~$A$ is noncompact
or else that~$n=2$.  Then there exists an embedded,
connected special Lagrangian submanifold~$L\subset\bbC^{n+1}$
that is $\SO(n)$-invariant and satisfies~$L\cap\bbC = A$.
\end{theorem}

\begin{proof}
Consider the set~$B$ that consists
of the pairs~$(z,P_\psi)$ with~$z\in A$ 
such that~$P_\psi$ is a special Lagrangian plane 
of the form~\eqref{eq: Ppsidef} that contains the line~$T_zA$.
By the discussion in~\S\ref{ssec: invsLagplanes}, 
it follows that~$B$ is a~$\bbZ_n$-bundle over~$A$.

If~$A$ is noncompact, then this bundle is trivial
and there exists a continuous 
(in fact, real-analytic)~$\psi:A\to\bbR$ such 
that~$(z,P_{\psi(z)})$ is a section of~$B$
over all of~$A$.
If~$A$ is compact and~$n=2$,
then the argument given at the end of the proof
of~Proposition~\ref{prop: Anonsing} shows that, 
because the rotation number of~$A$ is~$\pm1$, 
the bundle~$B$ is trivial in this case as well, 
so that, again, there is a 
mapping~$\psi:A\to \bbR/(\pi\bbZ)$ so that
that~$(z,P_{\psi(z)})$ is a section of~$B$
over all of~$A$.

Now, Proposition~\ref{prop: phi_exist},
together with the assumption that~$A$ is embedded,
implies that every point~$z\in A$ has an open 
neighborhood~$U_z\subset\bbC^{n+1}$ that intersects~$A$
in a connected arc~$U_z\cap A$ and in which there exists
an embedded, nonsingular, $SO(n)$-invariant special
Lagrangian submanifold, say~$L_z\subset U_z$ such
that~$L_z\cap\bbC = A\cap U_z$.  
By Theorem~\ref{thm: nuniqueL},
this~$L_z$ can be made unique by requiring that
$T_z(L_z) = P_{\psi(z)}$.  

By the connectedness of~$U_z\cap A$, and
the continuity of~$\psi$, it follows
that $T_w(L_z) = P_{\psi(w)}= T_w(L_w)$ 
for any~$w\in U_z\cap A$.
Consequently,~$T_y(L_z)=T_y(L_w)$ 
for all~$y\in U_z\cap U_w\cap A$.  
It then follows from Theorem~\ref{thm: nuniqueL} 
that~$L_z$ and~$L_w$ are equal in some open 
neighborhood of the interval~$U_z\cap U_w\cap A$.

It is now not difficult to conclude that there
exists an open neighborhood~$U$ of~$A$ itself and
an embedded, nonsingular, $SO(n)$-invariant 
special Lagrangian submanifold~$L\subset\bbC^{n+1}$
that agrees with~$L_z$ on some open neighborhood
of~$z$ for each~$z\in A$.
\end{proof}

\bibliographystyle{hamsplain}

\providecommand{\bysame}{\leavevmode\hbox to3em{\hrulefill}\thinspace}

\end{document}